\documentclass{amsart}

\usepackage{graphicx}
\usepackage{tikz}
\usepackage{tikz-cd}
\usepackage[latin1]{inputenc}

\tikzstyle{startstop}=[rectangle,rounded corners,draw=black,fill=green!20]
\tikzstyle{process}=[rectangle,draw=black,fill=green!20]
\tikzstyle{arrow}=[thick,->,>=latex]

\newtheorem{thm}{Theorem}[section]
\newtheorem{cor}[thm]{Corollary}
\newtheorem{lem}[thm]{Lemma}
\newtheorem{prop}[thm]{Proposition}
\newtheorem{exa}[thm]{Example}
\theoremstyle{definition}
\newtheorem{dfn}[thm]{Definition}
\theoremstyle{remark}
\newtheorem{rem}[thm]{Remark}
\numberwithin{equation}{section}

\newcommand{\dd}{d\mspace{-2.5mu}l l\mspace{-2.5mu}l}
\newcommand{\dL}{d\mspace{-2.5mu}l \mathbb{L}}
\begin{document}

\title[]{Knots in $S_{g} \times S^{1}$, their essential diagrams and virtual knots}%
\author{Seongjeong Kim}

\address{Seongjeong Kim, Jilin university}%
\email{kimseongjeong@jlu.edu.cn}%


\begin{abstract}
In \cite{Kim} it is shown that knots in $S_{g} \times S^{1}$ can be presented by virtual diagrams with a decoration, so called, {\em double lines}. In this paper, we study the essential diagram for each knot in $S_{g} \times S^{1}$, which has the minimal number of double lines. We prove that virtual knot theory is embedded in the theory of knots in $S_{g}\times S^{1}$. In the same time, one can obtain knots in $S^{2}\times S^{1}$ from 2-component links $L = K\sqcup T$ where $T$ is a trivial knot. By using knots in $S^{2} \times S^{1}$, we study the minimality and separability of such classical links.

\end{abstract}

\maketitle

\section{Introduction}

One of generalizations of classical knot theory is {\em virtual knot theory}.
\begin{dfn}
{\em A virtual link} is an equivalence class of virtual link diagrams modulo generalized Reidemeister moves described in Fig~\ref{vir_moves}. 
\begin{figure}
\centering\includegraphics[width=200pt]{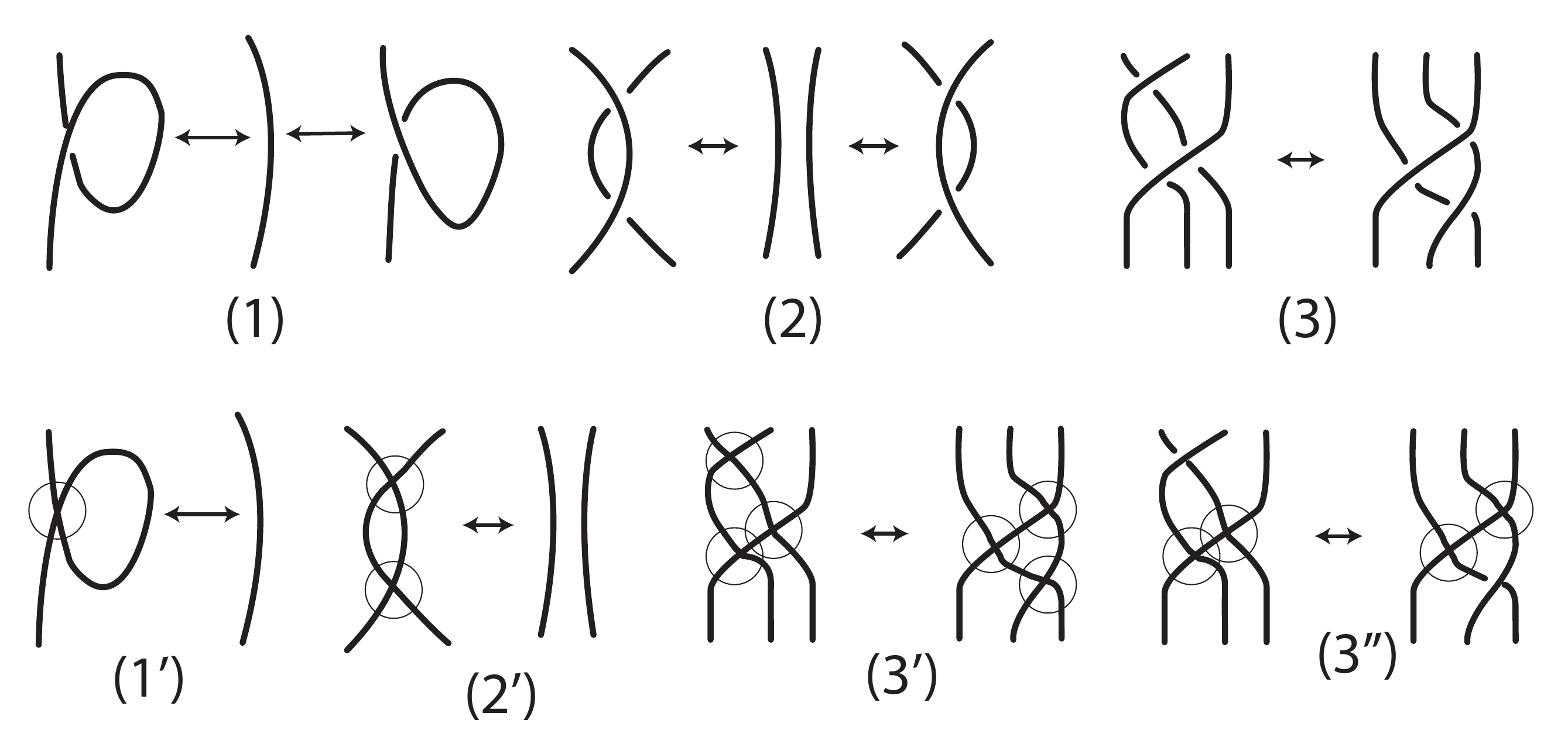} 
\caption{Generalized Reidemeister moves} \label{vir_moves}
 \end{figure}
\end{dfn}
It is well known that virtual links can be considered as links in a thickened surface $S_{g} \times [0,1]$ up to an equivalence relation, where $S_{g}$ is an oriented surface of genus $g$.

\begin{dfn}
{\em A virtual link} is a smooth embedding $L$ of a disjoint union of $S^{1}$ into $S_{g} \times [0,1]$. Each image of $S^{1}$ is called {\em a component} of $L$. A link of one component is called {\em a virtual knot}.

\end{dfn}

\begin{dfn}
  Let $L$ and $L'$ be two virtual links. If $L'$ can be obtained from $L$ by isotopy and stabilization/destabilization of $S_{g} \times [0,1]$, then we call $L$ and $L'$ are \textit{equivalent}.
\end{dfn}

It is well known that {\bf those two definitions of the virtual knot are equivalent,} for example, see \cite{CarterKamadaSaito, KamadaKamada, Kauffman}.

In virtual knot theory, by using {\em the parity} defined by V.O. Manturov, many invariants for classical knots can be non-trivially extended to invariants for virtual knots. It gives several interesting geometrical properties, for details, see \cite{Manturov_vir_book}. 

One of important properties, shown by using parity, is that classical knot theory is embedded in virtual knots. More precisely, it states that if two classical knot diagrams $D_{1}$ and $D_{2}$ are equivalent modulo generalized Reidemeister moves, then they are equivalent modulo Reidemeister moves. Sometimes we simply say that {\em classical knots are embedded in virtual knots.}

 In \cite{CM}, M. Chrisman and V.O. Manturov studied virtual knots by using 2-component classical link $K \sqcup J$ in $S^{3}$ with $lk(K,J)=0$, where $J$ is a fibered knot. They obtain a virtual knot $\hat{K}$ from $K \sqcup J$ with $lk(K,J)=0$ as a lifting of $K \subset S^{3}\backslash J \cong \Sigma_{J} \times [0,1]/x\sim f(x)$ to $\Sigma_{J} \times \mathbb{R}$, where $f :\Sigma_{J} \rightarrow \Sigma_{J}$ is a homeomorphism. 
 
 We are interested in links in $S_{g}\times S^{1}$, where $S_{g}$ is an oriented surface of genus $g$. In \cite{Kim}, it is proved that links in $S_{g} \times S^{1}$ can be presented by virtual knot diagrams with decorations, called {\em double lines}. Moreover the isotopy of two knots in $S_{g} \times S^{1}$ can be presented by finitely many local moves. It is an extension of the work in \cite{DabkowskiMroczkowski}.  In \cite{Kim-winding-parity}, the author defined the winding parity, which is a similar notion to parity in the virtual knot theory for knots in $S_{g}\times S^{1}$ and applied to distinguish some families of knots in $S_{g}\times S^{1}$.

One can notice that there is a natural map from virtual knots (knots in $S_{g} \times [0,1]$) to knots in $S_{g} \times S^{1}$.
This paper is contributed to prove that virtual knot theory is embedded in the theory of knots in $S_{g} \times S^{1}$. In Section 2, we introduce the basic notions of links in $S_{g} \times S^{1}$. In Section 3, we show that the virtual knot theory is embedded in the theory of knots in $S_{g} \times S^{1}$. In Section 4, we extend the idea in Section 3 for knots in degree $k$ where $k \neq 0$. We find a diagram with the minimal number of double lines equivalent to the original diagram, which is called {\em an essential diagram}. In Section 5, we apply statements proved in Section 3 and Section 4. First, we expand our classification of knots in $S^{2} \times S^{1}$, which is studied in \cite{Kim-winding-parity}. And then we describe a way to obtain diagrams with double lines from classical links of 2-components with a trivial knot component. It provides a sufficient condition for 2-component links with a trivial component to be separable. 

\section{Link in $S_{g} \times S^{1}$ and its diagrams}

Let $S_{g}$ be an orientable surface of genus $g$. Let us define links in $S_{g} \times S^{1}$ analogously to virtual links by using underlying surfaces as follows:
\begin{dfn}
Let $S_{g}$ be an oriented surface of genus $g$. {\em A link $L$  in $S_{g} \times S^{1}$} is a pair of $S_{g} \times S^{1}$ and a smooth embedding of a disjoint union of $S^{1}$'s into $S_{g} \times S^{1}$. We denote it by $(L, S_{g} \times S^{1})$. Each image of $S^{1}$ in $S_{g} \times S^{1}$ is called {\em a component} of $L$. A link of one component is called {\em a knot in $S_{g} \times S^{1}$.}
\end{dfn}

\begin{dfn}
Let $(L, S_{g} \times S^{1})$ and $(L', S_{g'} \times S^{1})$ be two links in $S_{g} \times S^{1}$ and $S_{g'} \times S^{1}$. We call $(L, S_{g} \times S^{1})$ and $(L', S_{g'} \times S^{1})$ are {\em equivalent,} if $(L', S_{g'} \times S^{1})$ can be obtained from $(L, S_{g} \times S^{1})$ by isotopy and stabilization/destabilization of $S_{g} \times S^{1}$.
\end{dfn}
By the {\em destabilization for $(L, S_{g} \times S^{1})$} we mean the following:
let $C$ be a non-contractible circle on the surface $S_{g}$ such that there exists a torus $T$ in $S_{g} \times S^{1}$ homotopic to the torus $C \times S^{1}$ and not intersecting the link. Then the destabilization is cutting of the manifold $S_{g} \times S^{1}$ along the torus $C \times S^{1}$ and pasting of two newborn components of boundary by $D^{2} \times S^{1}$.
The {\em stabilization for $S_{g} \times S^{1}$} is the converse operation to the destabilization.

First let us construct diagrams for $(L, S_{g}\times S^{1})$ on the surface $S_{g}$ as follows: let $L$ be an (oriented) link in $S_{g} \times S^{1}$. Assume that an orientation is given on $S^{1}$. Suppose that $x_{0} \in S^{1}$ is a point such that $S_{g} \times \{x_{0}\} \cap L$ is a set of finite points with no transversal points. 
\begin{figure}[h]
\begin{center}
 \includegraphics[width = 12cm]{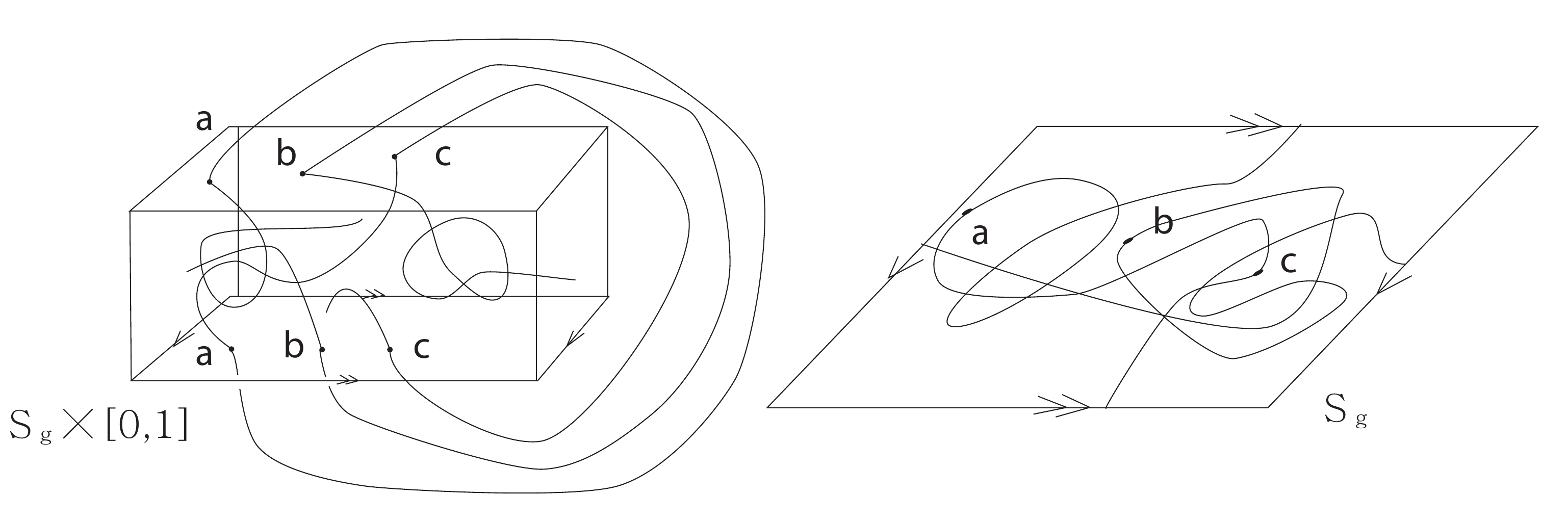}

\end{center}

\caption{Schematic figures of links in $S_{g} \times S^{1}$ and its projection on the surface $S_{g}\times \{0\}$}\label{fig:SS-diag-rel-1}
\end{figure}
Then there exists a natural diffeomorphism $f$\label{def:cutting-map} from $(S_{g} \times S^{1} )\backslash$  $(S_{g} \times \{x_{0}\})$ to $S_{g} \times (0,1) \subset S_{g} \times [0,1]$. Let $M_{L}$ $=$ $\overline{f((S_{g} \times S^{1}) - (S_{g} \times \{x_{0}\}))} \cong S_{g} \times [0,1]$. Then $\overline{f(L)}$ in $M_{L}$ consists of finitely many circles and arcs with exactly two boundaries on $S_{g} \times \{0\}$ and $S_{g} \times \{1\}$.

Let $D_{\overline{f(L)}}$ be the image of a projection of $\overline{f(L)}$ on the $S_{g} \times \{0\}$. Notice that boundaries of $\overline{f(L)}$ are paired to be projected to the same point on $S_{g}\times \{0\}$, we call them {\em vertices}. It follows that the diagram $D_{\overline{f(L)}}$ of $L$ on $S_{g}$ has vertices corresponding to two boundary points on $S_{g} \times \{0\}$ and $S_{g} \times \{1\}$ of $\overline{f(L)}$ as described in the right of Fig.~\ref{fig:SS-diag-rel-1}.

\begin{figure}[h!]
\begin{center}
 \includegraphics[width = 8cm]{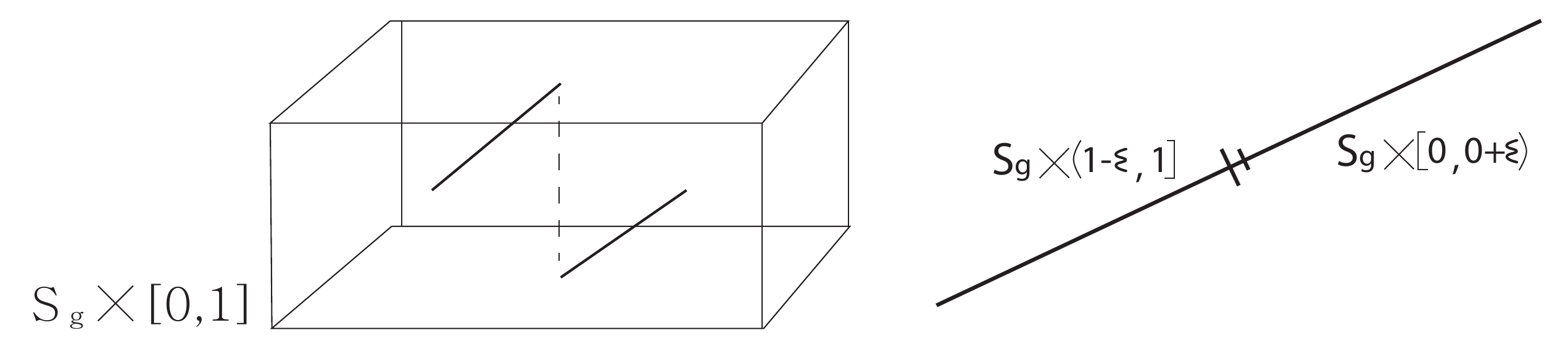}

\end{center}
\caption{Local image near to an intersection of a link with $S_{g} \times \{x_{0}\}$ and a corresponding double line}\label{fig:Vertex}
\end{figure}

Notice that two arcs near to a vertex are the images of arcs near to $S_{g} \times \{0\}$ and $S_{g} \times \{1\}$ in $M_{L} \cong S_{g} \times [0,1]$, respectively, as described in Fig.~\ref{fig:Vertex}. We change the point to two small lines so that if one of the lines is connected with an arc which is near to $S_{g} \times \{1\}$, then the line is longer than the another, as describe in Fig.~\ref{fig:Vertex}.

Since $D_{\overline{f(L)}}$ is a framed 4-valent graph with double lines on the surface $S_{g}$, which comes from $\overline{f(L)}$ in $S_{g} \times [0,1]$, one can give classical crossing information for each 4-valent vertex. That is, a link $L$ in $M_{L}$ has a knot diagram with double lines on $S_{g}$. Simply we call it {\em a diagram on $S_{g}$ with double lines.} 

\begin{prop}[M.K. Dabkowski, M. Mroczkowski (2009) \cite{DabkowskiMroczkowski}, S. Kim (2018) \cite{Kim}]\label{thm:diag_on_surface}
   Let $(L, S_{g} \times S^{1})$ and $(L', S_{g} \times S^{1})$ be two links in $S_{g} \times S^{1}$. Let $D_{L}$ and $D_{L'}$ be diagrams of $L$ and $L'$ on $S_{g}$, respectively. Then $(L, S_{g} \times S^{1})$ and $(L', S_{g} \times S^{1})$ are isotopic if and only if $D_{L'}$ can be obtained from $D_{L}$ by applying finitely many moves in Fig.~\ref{fig:moves1}.

  \begin{figure}[h!]
\begin{center}
 \includegraphics[width = 9cm]{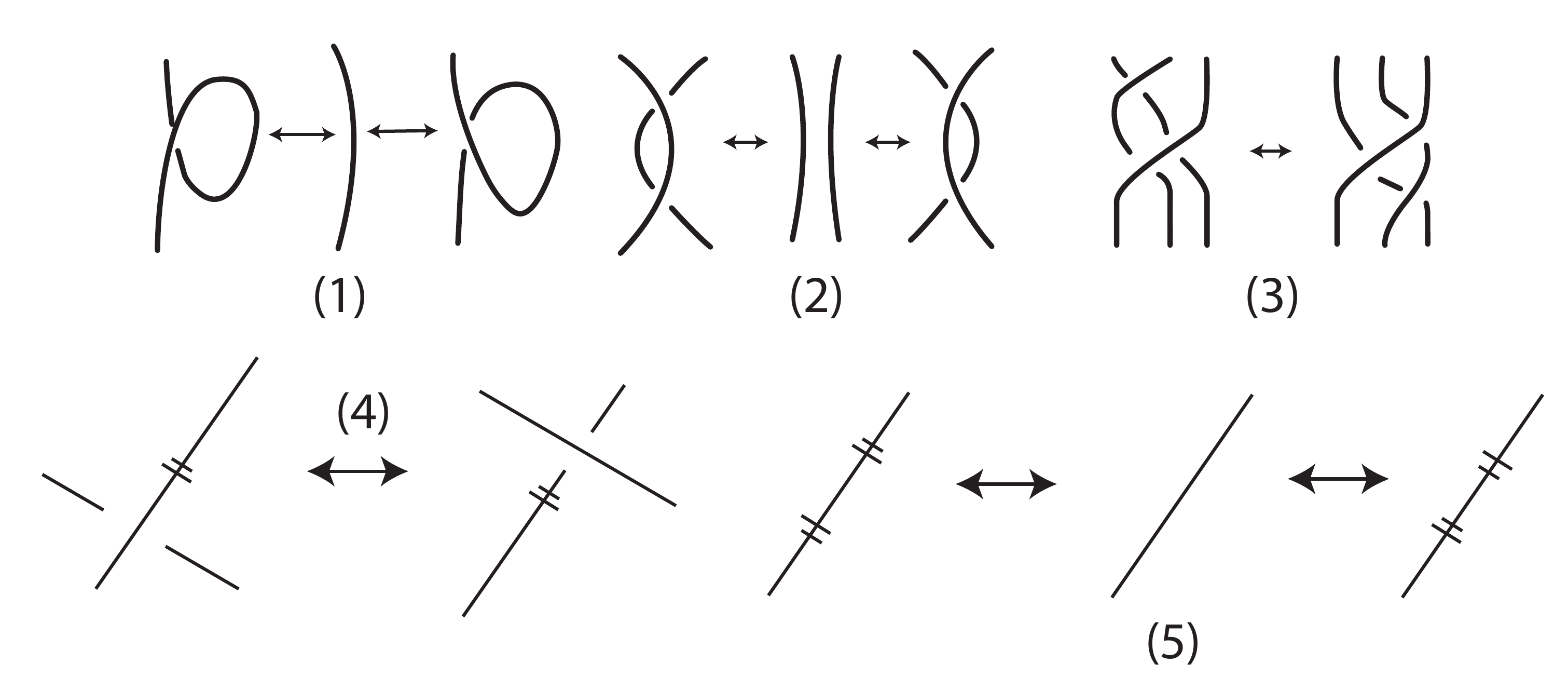}

\end{center}

 \caption{Moves for diagrams on $S_{g}$}\label{fig:moves1}
\end{figure}
\end{prop}

\begin{cor}
Let $(L, S_{g} \times S^{1})$ and $(L', S_{g-1} \times S^{1})$ be two links in $S_{g} \times S^{1}$ and $S_{g-1} \times S^{1}$, respectively. Let $D_{L}$ and $D_{L'}$ be diagrams $S_{g}$ and $S_{g-1}$ of $L$ and $L'$, respectively. Then $(L', S_{g-1} \times S^{1})$ is obtained from $(L, S_{g} \times S^{1})$ by a destabilization if and only if $(D_{L'}, S_{g-1})$ can be obtained from $(D_{L}, S_{g})$ by a destabilization of $S_{g}$.
\end{cor}

\begin{figure}[h]
\begin{center}
 \includegraphics[width = 12cm]{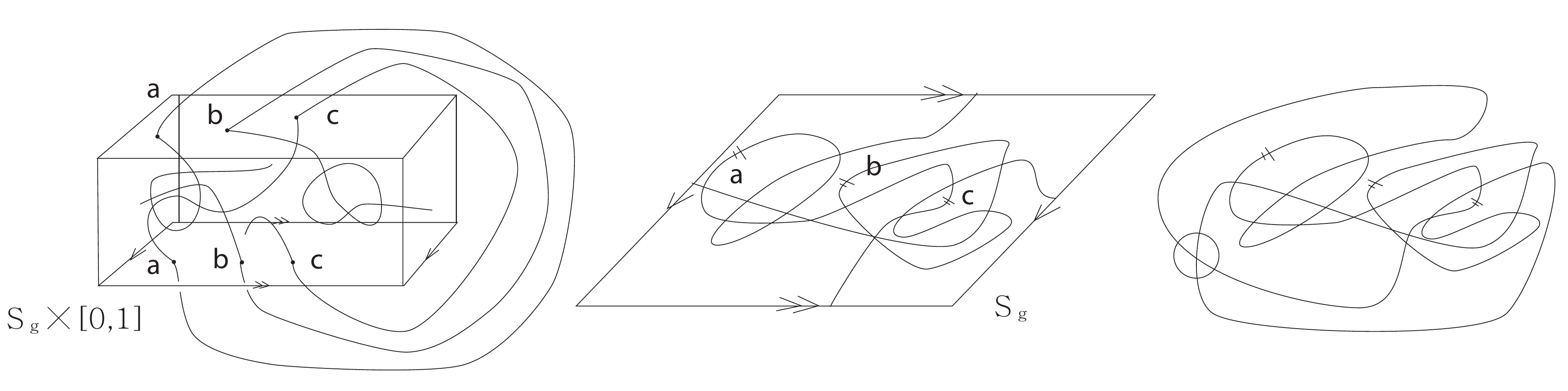}

\end{center}

\caption{Schematic figures of links in $S_{g} \times S^{1}$ and its projection on the plane}\label{fig:SS-diag-rel-2}
\end{figure}

Now, let us construct diagrams for links in $S_{g}\times S^{1}$ on the plane by using diagrams on $S_{g}$. For a link in $S_{g} \times S^{1}$ let a diagram $D$ on $S_{g}$ with double lines be given. We may assume that the diagram is drawn on $2g$-gon presentation of $S_{g}$ as in the middle of Fig.~\ref{fig:SS-diag-rel-2}. Connect points on boundaries of $2g$-gon with corresponding points by arcs outside $2g$-gon. By changing intersections between arcs outside $2g$-gon to virtual crossings, we obtain a diagram with double lines and virtual crossings, see the right in Fig.~\ref{fig:SS-diag-rel-2}. We call it {\it a diagram with double lines on the plane} for links in $S_{g}\times S^{1}$, or simply {\it a diagram for links in $S_{g}\times S^{1}$}.

The following theorem holds.
\begin{prop}[Kim (2018) \cite{Kim}]\label{thm:diag_on_plane}
   Let $(L, S_{g} \times S^{1})$ and $(L', S_{g'} \times S^{1})$ be two links. Let $D_{L}$ and $D_{L'}$ be diagrams of $L$ and $L'$ on the plane, respectively. Then $L$ and $L'$ are equivalent if and only if $D_{L'}$ can be obtained from $D_{L}$ by applying finitely many moves in Fig.~\ref{fig:moves2}.

  \begin{figure}[h!]
\begin{center}
 \includegraphics[width = 12cm]{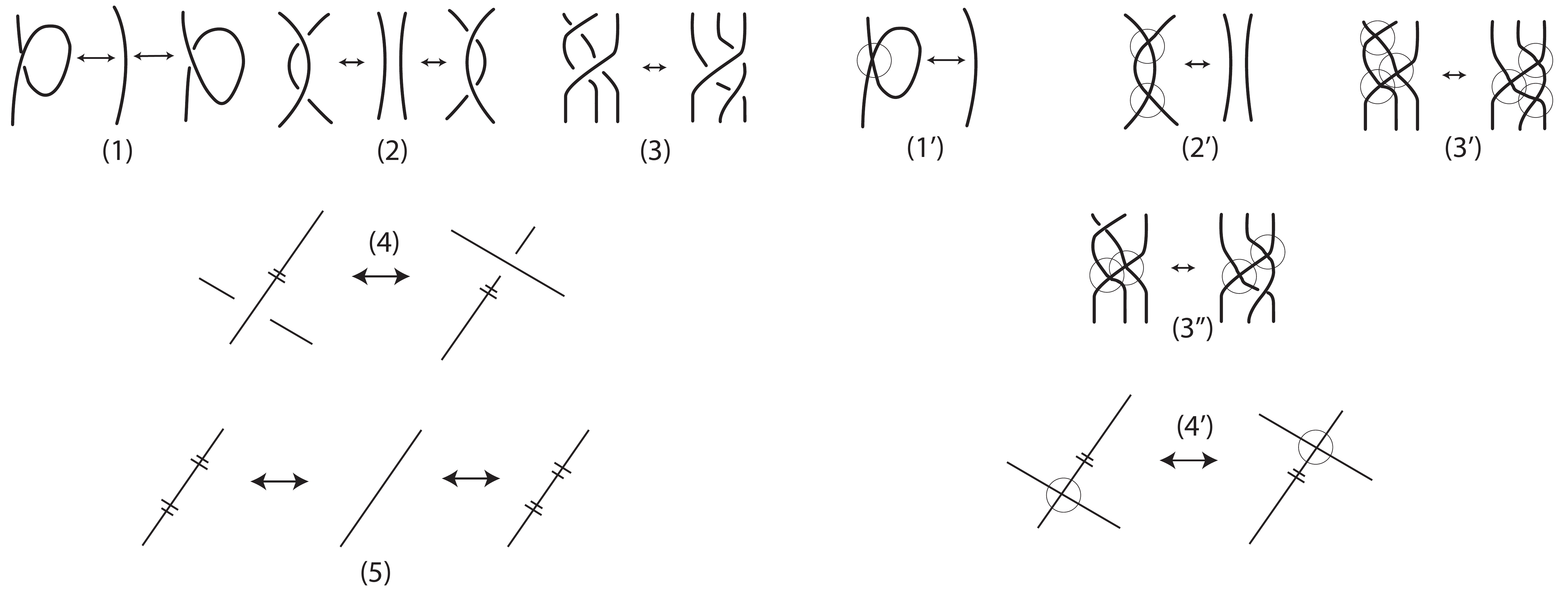}

\end{center}

 \caption{Moves for links in $S_{g}\times S^{1}$}\label{fig:moves2}
\end{figure}
\end{prop}

\begin{figure}[h]
\begin{center}
 \includegraphics[width = 9cm]{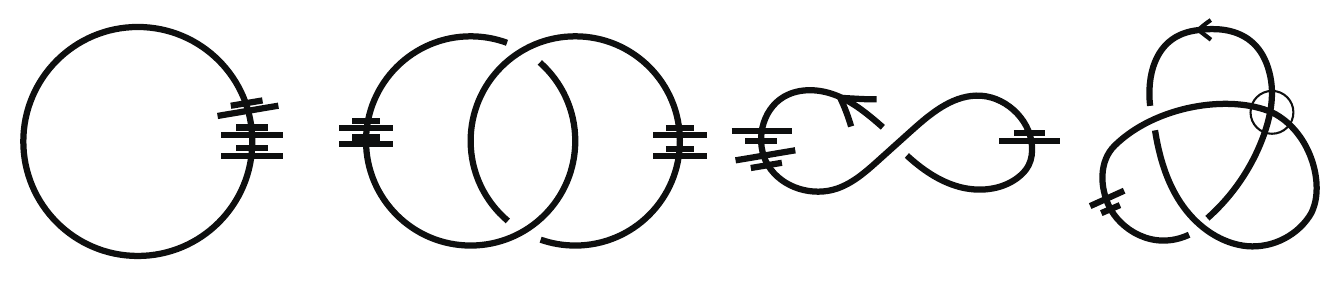}

\end{center}

\caption{Diagrams with double lines}\label{fig:exa-diag-dl}
\end{figure}

Simply speaking, knots and links in $S_{g}\times S^{1}$ can be presented by virtual knot diagrams with a decoration, which is the {\em double line} as described in Fig.~\ref{fig:exa-diag-dl}.

For a diagram $D_{L}$ with double lines, we call the virtual diagram, where double lines appear for $D$, {\em the base diagram}, and $D_{L}$ is called {\em to be based on $D$.}\\[2mm]

{\bf Terminology 1.} In this paper, arcs of a diagram with double lines between two double lines are called {\em long arcs}, arcs between crossings are called {\em arcs} and arcs between a crossing and a double lines are called {\em short arcs}.\\[2mm]

{\bf Terminology 2.} By the move (5) in Fig. \ref{fig:moves2} two neighboring double lines appear or disappear as a pair. We call this move {\em double line appearing} and {\em double line canceling}.

\begin{rem}\label{rem:cro-change}
As described in Fig.~\ref{cro_change}, by adding two double lines one can change over/under information of a crossing. In the present paper, we call such a move {\em a crossing change move}. Note that the move (4) in Fig.~\ref{fig:moves2} can be obtained by using move (5) and the crossing change move. Sometimes, we replace the move (4) by crossing change move.

  \begin{figure}[h!]
\begin{center}
 \includegraphics[width = 8cm]{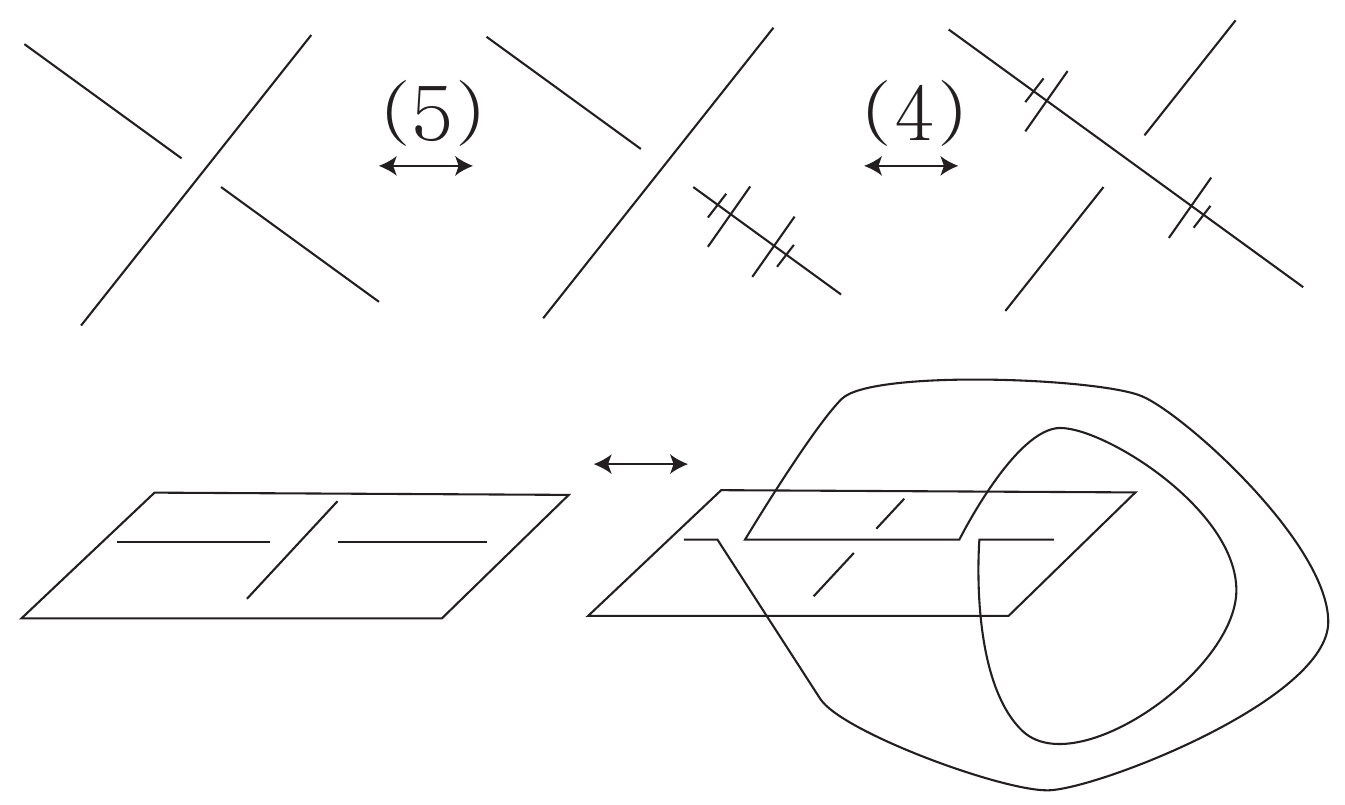}

\end{center}

 \caption{A crossing change with additional two double lines with additional two double lines}\label{cro_change}
\end{figure}
\end{rem}


\subsection{Degree of knots in $S_{g}\times S^{1}$} 
From now on we are mainly interested in {\it oriented knots} in $S_{g} \times S^{1}$.  Let $\Pi: \mathbb{R} \rightarrow S^{1}$ be the covering defined by $\Pi(r) = e^{2\pi r i}$. Then the function $Id_{S_{g}} \times \Pi : S_{g}\times \mathbb{R} \rightarrow S_{g}\times S^{1}$ is also a covering over $S_{g}\times S^{1}$, where $Id_{S_{g}} : S_{g} \rightarrow S_{g}$ is the identity map.

\begin{center}
\begin{tikzcd} 
&|[alias=Z]| S_{g}\times \mathbb{R} \arrow[r,"\phi_{2}"]\arrow[d, "Id_{S_{g}}\times \Pi"]& \mathbb{R} \arrow[d,"\Pi"]\\
S^{1}\arrow[to=Z, "\hat{K}"]\arrow[r, "K"] &S_{g}\times S^{1} \arrow{r}&S^{1}
\end{tikzcd}
\end{center}

Let $K : [0,1] \rightarrow S_{g}\times S^{1}$ be a knot with $K(0)=K(1)$.
 Let $\tilde{K}$ be a lifting of $K$ into $S_{g} \times \mathbb{R}$ along the covering $Id_{S_{g}} \times \Pi : S_{g}\times \mathbb{R} \rightarrow S_{g}\times S^{1}$. When $\phi_{2} \circ \hat{K}(0) =0$ for the projection $\phi_{2} : S_{g} \times \mathbb{R} \rightarrow \mathbb{R}$, {\em the degree $deg(K)$ of a knot $K$ in $S_{g} \times S^{1}$} is defined by 
 $$deg(K) = \phi_{2} \circ \hat{K}(1)\in \mathbb{Z}.$$
 It is easy to see that the degree $deg(K)$ of a knot $K$ in $S_{g} \times S^{1}$ is an invariant for knots in $S_{g} \times S^{1}$.

\begin{figure}[h!]
\begin{center}
 \includegraphics[width = 5cm]{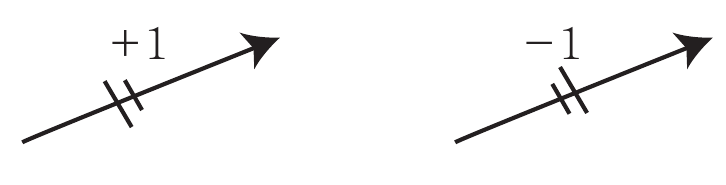}
\caption{Signs of double lines}\label{fig:signs-doublelines}
\end{center}
\end{figure}

\begin{figure}[h!]
\begin{center}
 \includegraphics[width = 8cm]{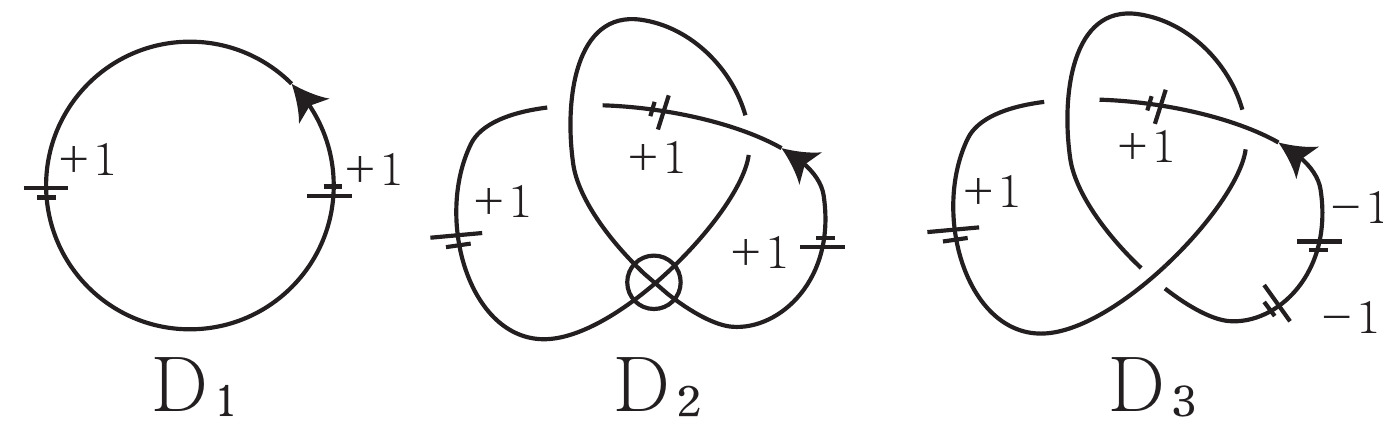}
\caption{$D_{1}$, $D_{2}$ and $D_{3}$ are oriented diagrams with double lines. Each double line has a sign. We obtain that $deg(D_{1}) = 2$, $deg(D_{2}) = 3$, $deg(D_{3}) = 0$.}\label{fig:exa-diag-degree}
\end{center}
\end{figure}

 The degree of a given oriented knot $K$ in $S_{g} \times S^{1}$ can be calculated by using a diagram with double lines in the following way:
Let $D$ be an oriented diagram with double lines of an oriented knot $K$ in $S_{g} \times S^{1}$. Let us give $\pm 1$ to double lines with respect to the orientation as described in Fig.~\ref{fig:signs-doublelines}, call it {\em a sign of a double line}. 
Then the degree of $K$ is equal to the sum of signs of all double lines, for example, see Fig.~\ref{fig:exa-diag-degree}.

\subsection{Winding parity for crossings of a knot in $S_{g} \times S^{1}$}\label{sec:label}
\begin{figure}[h!]
\begin{center}
 \includegraphics[width = 3.5cm]{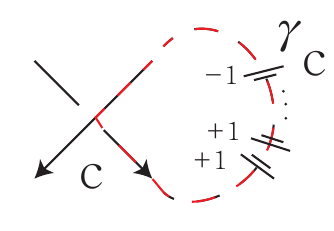}

\end{center}
 \caption{A winding of a crossing}\label{fig:label-crossing}
\end{figure}

For a crossing $c$ take a closed curve $\gamma_{c}$ starting from under-crossing of $c$, see Fig.~\ref{fig:label-crossing}. {\em A winding parity of $c$} valued in $\mathbb{Z}_{deg(K)}$ is the sum of $\pm 1$'s for double lines on $\gamma_{c}$.
\begin{figure}[h!]
\begin{center}
 \includegraphics[width = 6cm]{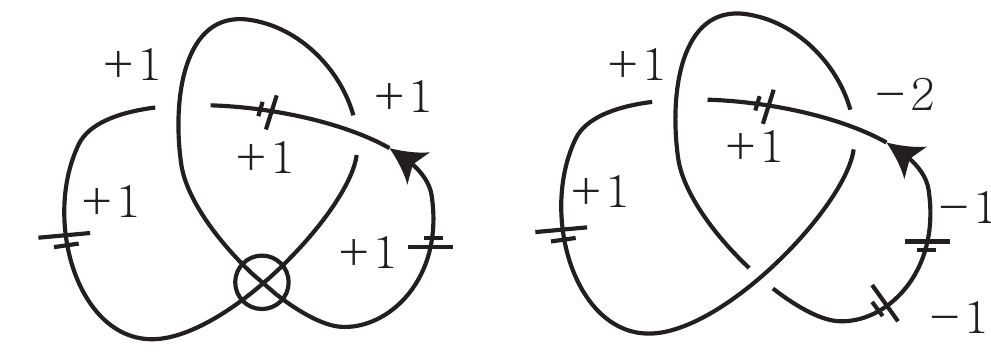}
 \put(-132, -10){$D_{1}$}
 \put(-45, -10){$D_{2}$}
 \caption{Examples of winding parities of crossings of diagrams with double lines. Note that, since $deg(D_{1})=3$ and $deg(D_{2})=0$, winding parities for crossings of $D_{1}$ are valued in $\mathbb{Z}_{3}$, but are valued in $\mathbb{Z}$ for crossing of $D_{2}$.}
\end{center}
\end{figure}

\begin{lem}
 The winding parity for crossings satisfy the properties described in Fig.~\ref{fig:labelmoves2}.\label{lem:labelmoves}

  \begin{figure}[h!]
\begin{center}
 \includegraphics[width = 8.5cm]{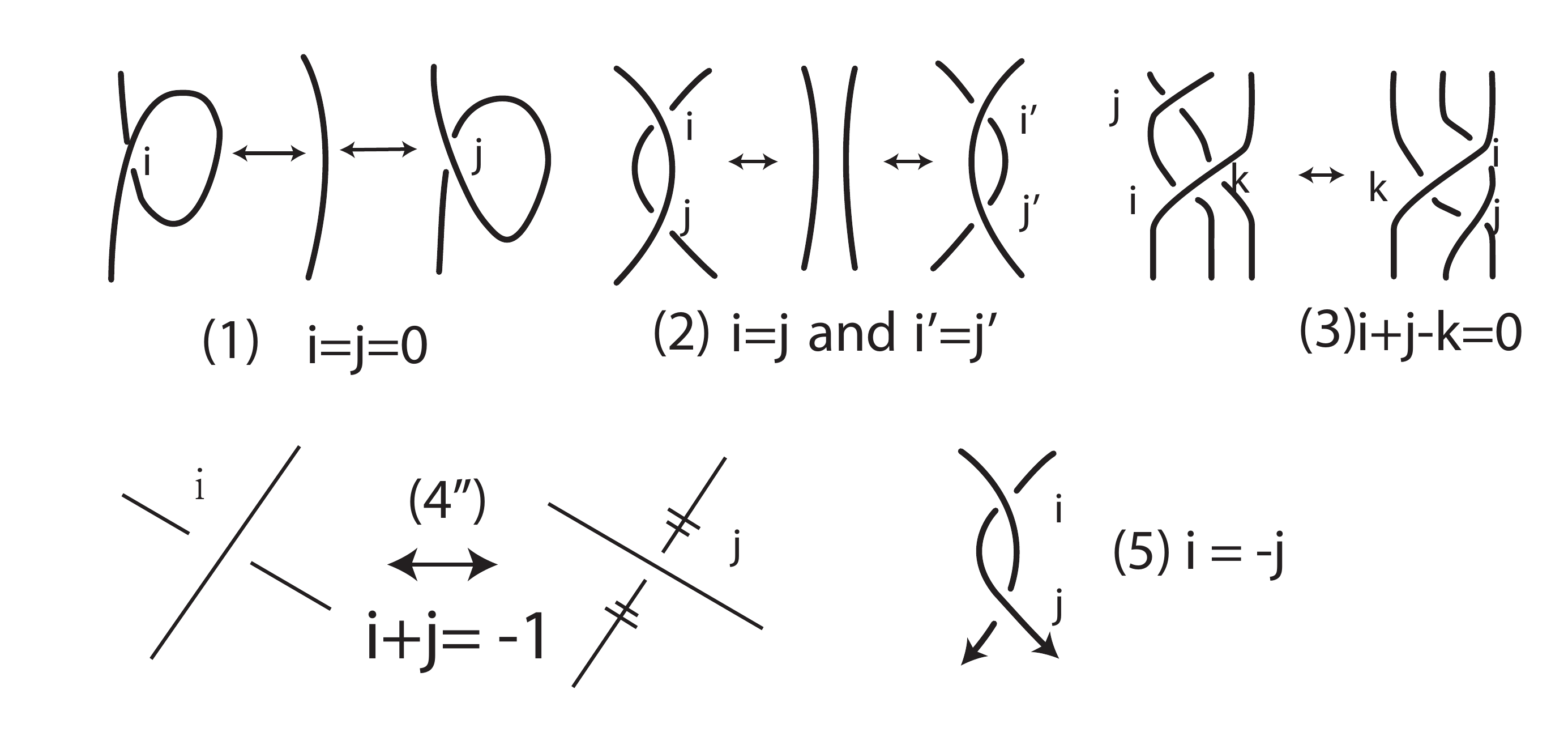}

\end{center}
 \caption{Properties of labels of crossings}\label{fig:labelmoves2}
\end{figure}
\end{lem}

\begin{rem}
Geometrically, the winding parity of a crossing $c$ means how many times the curve from the crossing $c$ to itself turns around $S^{1}$.
\end{rem}

\begin{rem}
    In \cite{Kim-winding-parity}, the winding parity is defined for more general notions with conditions derived from the properties in Lemma~\ref{lem:labelmoves}, but the present one introduced above is called {\bf a label}, which satisfies the conditions for winding parity. In this paper, the winding parity means the one defined above.
\end{rem}






\section{Embedding of virtual knots into knots in $S_{g}\times S^{1}$}
In this section we show that virtual knot theory is embedded well into the theory of knots in $S_{g} \times S^{1}$ in the following sense: two virtual knot diagrams are equivalent modulo moves for knots in $S_{g} \times S^{1}$, then they are equivalent modulo generalized Reidemeister moves. To achieve the above goal, we introduce one local move for knots in $S_{g}\times S^{1}$. 

\begin{lem}
    From the moves for diagrams with double lines we obtain the crossing sliding move described in Fig.~\ref{fig:cro-sliding}.
\begin{figure}
    \centering
    \includegraphics[width=0.5\linewidth]{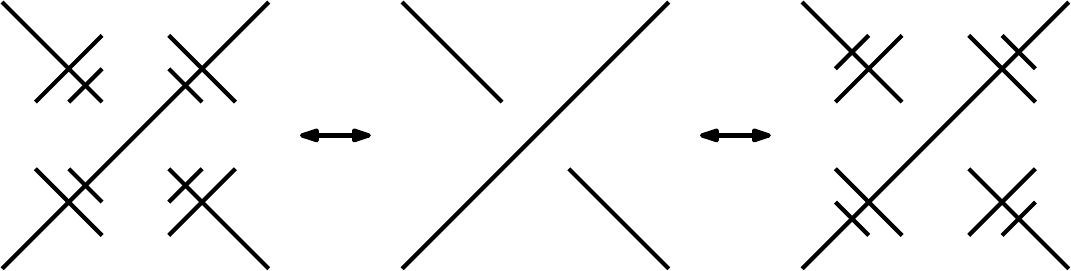}
    \caption{Crossing sliding move}
    \label{fig:cro-sliding}
\end{figure}
\end{lem}

\begin{proof}
    The crossing sliding move is obtained by applying crossing change move twice.
\end{proof}

\begin{figure}[h]
\begin{center}
\includegraphics[width = 3cm]{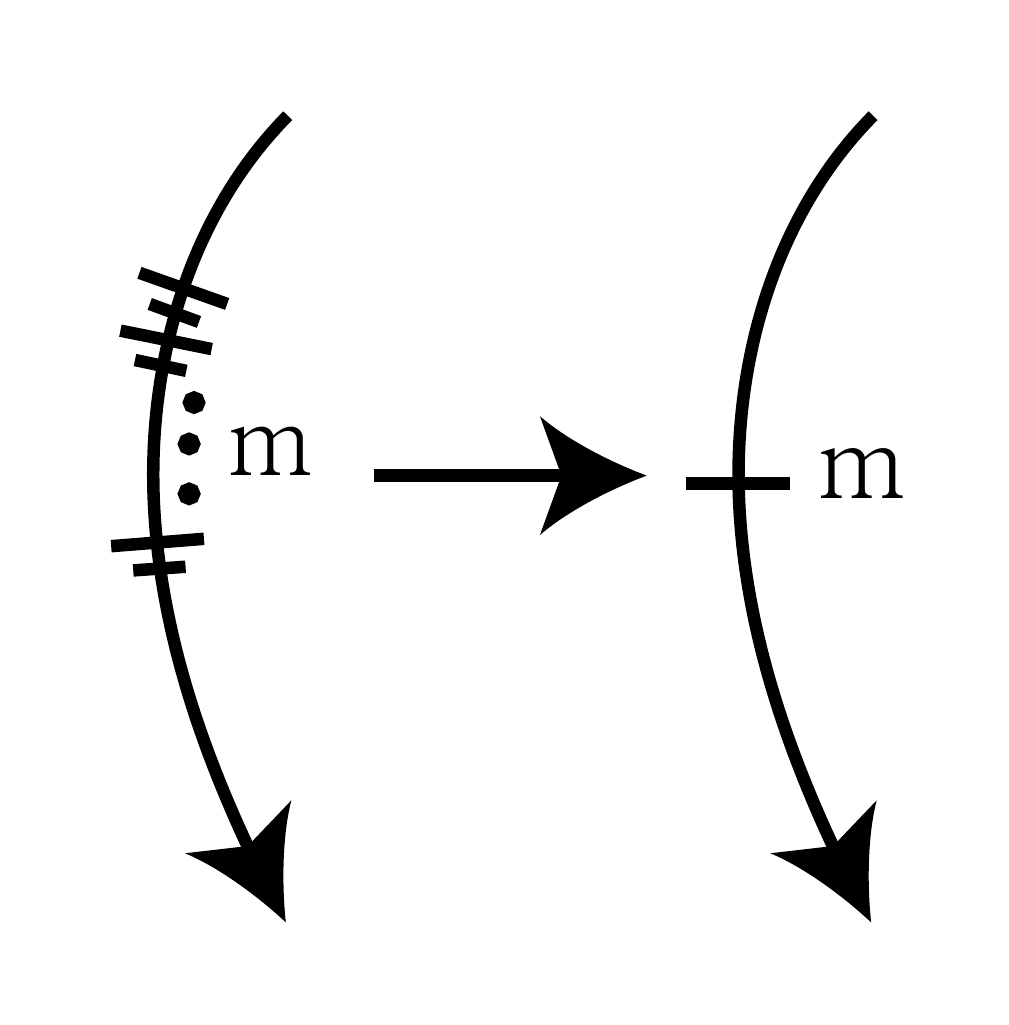}
\end{center}
\caption{The neighboring double lines with sum of signs of neighboring double lines on an arc is presented by one line with the integer number $m$}\label{fig:label-doubleline-arc}
\end{figure}

For convenience, we present neighboring double lines by one line with an integer number $m$, where $m$ is the sum of signs of neighboring double lines, as illustrated in Fig. \ref{fig:label-doubleline-arc}. 

Let $\mathcal{D}_{dl}$ be the set of equivalence classes of oriented diagrams with double lines. Let $\mathcal{D}_{dl}^{k}$ be the set of equivalence classes of diagrams with double lines with degree $k$ for $k \in \mathbb{Z}$. By definition, it is clear that $\mathcal{D}_{dl} = \sqcup \mathcal{D}_{dl}^{k}$.

\begin{figure}
    \centering
    \includegraphics[width=0.4\linewidth]{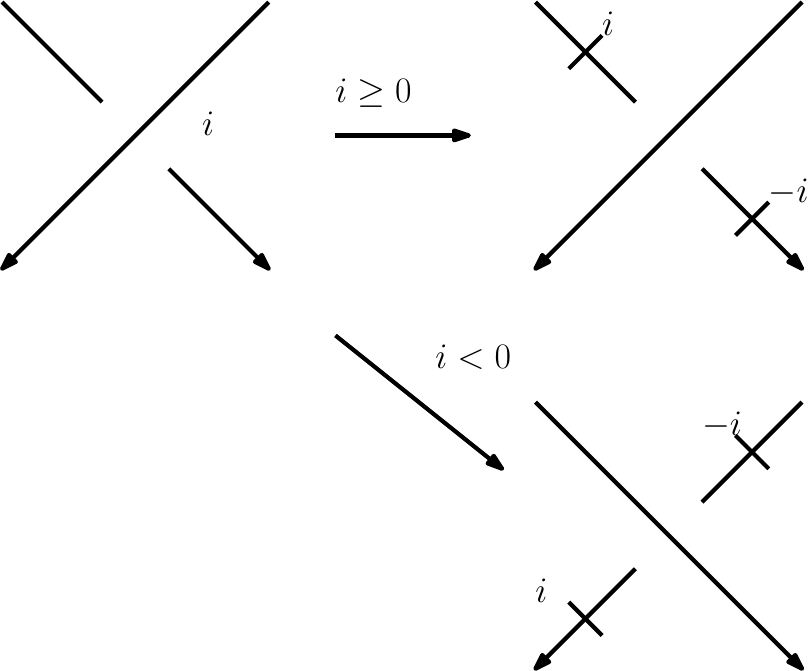}
    \caption{The winding parity projection}
    \label{fig:wp-proj}
\end{figure}

Let us define a map from $\mathcal{D}_{dl}^{0}$ to $\mathcal{D}_{dl}^{0}$ as follows: for each crossing $c$ with winding parity $i$ of an oriented diagram $D$ in $\mathcal{D}_{dl}^{0}$, if $i\geq 0$, then we add $i$ negative double lines and $i$ positive double lines to outgoing and incoming under arcs respectively. If $i<0$, then we change over/under crossing and add $i$ positive double lines and and $i$ negative double lines to outgoing and incoming under arcs respectively. See Fig.~\ref{fig:wp-proj}. 
Let us denote the obtained diagram by $pr_{wp}(D)$. It is clear that the degree of $pr_{wp}(D)$ is $0$, that is, $pr_{wp}(D) \in \mathcal{D}_{dl}^{0}$. 

\begin{thm}\label{thm:wp-projection}
    Let $pr_{wp} : \mathcal{D}_{dl}^{0} \rightarrow \mathcal{D}_{dl}^{0}$ be the map defined as above. Then $pr_{wp}$ is well-defined. In particular, if $D_{1}$ and $D_{2}$ are equivalent as knots in $S_{g} \times S^{1}$, then $pr_{wp}(D_{1})$ and $pr_{wp}(D_{2})$ are equivalent up to the crossing sliding move and the moves in Fig. \ref{fig:moves2} except for the move (4).
\end{thm}

\begin{figure}
    \centering
    \includegraphics[width=0.4\linewidth]{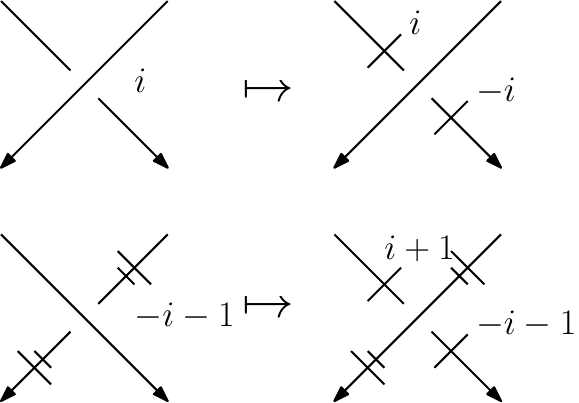}
    \caption{The equivalence relation up to the crossing change move is preserved by the winding parity projection. Here we use the crossing sliding move and the double line canceling.}
    \label{fig:pf-move4}
\end{figure}

\begin{proof}
    Let $D$ and $D'$ be two diagrams with double lines of degree $0$. Suppose that $D'$ is obtained from $D$ by one of moves for knots in $S_{g} \times S^{1}$. If $D'$ is obtained from $D$ by the first Reidemeister moves, double lines canceling moves or moves containing virtual crossings, then it is easy to see that $pr_{wp}(D')$ can be obtained from $pr_{wp}(D)$ by moves for knots in $S_{g}\times S^{1}$ except for the crossing change move. Suppose that $D'$ is obtained from $D$ by the second Reidemeister move. Then two crossings have the same winding parity $i$. If $i\geq 0$, then we remain these crossings, but add double lines on incoming and outgoing under-arcs for each crossing as described in Fig.~\ref{fig:pf-move2}.
    \begin{figure}
    \centering
\includegraphics[width=0.35\linewidth]{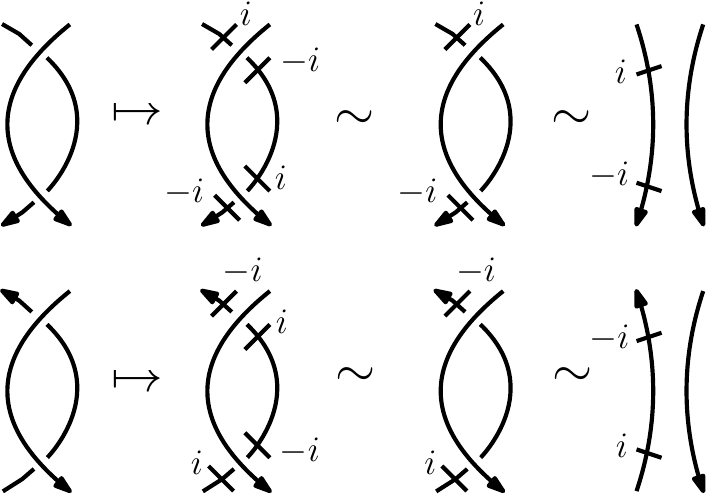}
    \caption{The equivalence relation up to the move (2) is preserved by the winding parity projection. Note that we use the double line canceling.}
    \label{fig:pf-move2}
\end{figure}
Then by applying double line cancelings and the second Reidermeister moves, $pr_{wp}(D')$ is obtained from $pr_{wp}(D)$. 
Suppose that $D'$ is obtained from $D$ by the third Reidemeister move. Say, three crossings in the move have winding parity $i,j,k$ as described in Fig.~\ref{fig:pf-move3}. Note that they satisfy the equality $i+j-k=0$.
\begin{figure}
    \centering
    \includegraphics[width=0.5\linewidth]{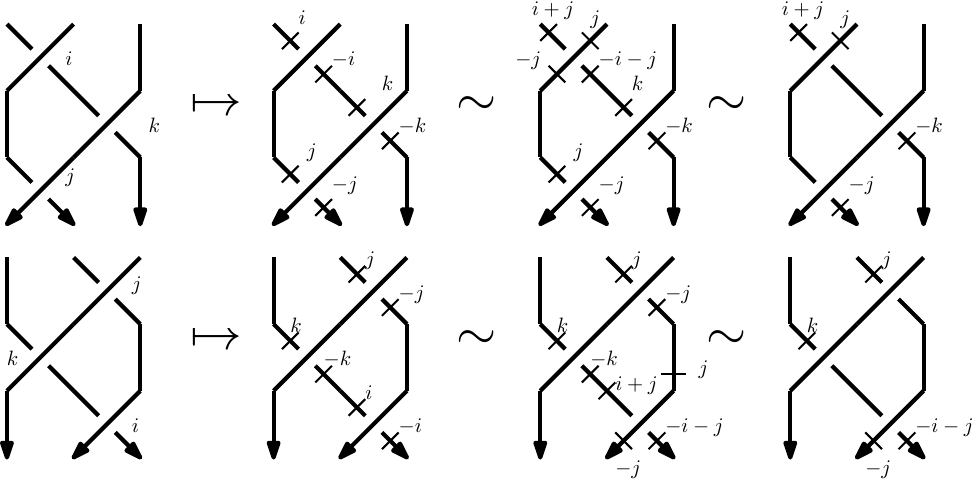}
    \caption{The equivalence relation up to the move (3) is preserved by the winding parity projection. Here we use the crossing sliding move and the double line canceling. In particular, we use the property of winding parities that $i+j-k=0$}
    \label{fig:pf-move3}
\end{figure}
We have six possible cases with respect to the signs of $i,j,k$. If $i,j,k \geq 0$, then we remain crossings, but add double lines on under-arcs with respect to $i,j,k$. By applying crossing sliding move $j$ times as described in Fig.~\ref{fig:pf-move3}. Since $i+j-k = 0$, or equivalently, $i+j=k$, we can remove all double lines on the triangle, where the third Reidemeister move is applied, by double line canceling. 
Finally, for the move (4), we replace it by the crossing change move. Suppose that $D'$ is obtained from $D$ by the crossing change move. Say, the diagram $D$ has a crossing $c$ with a winding parity $i \geq 0$, but $D'$ has a crossing $c'$, which is obtained from $c$ by crossing change move, with the winding parity $-i-1$, which is negative. In Fig.~\ref{fig:pf-move4} the upper right diagram describes $pr_{wp}(D)$ and the lower right diagram is  $pr_{wp}(D')$. Then $pr_{wp}(D)$ is obtained from $pr_{wp}(D)$ by applying crossing sliding move once. 
\end{proof}

Note that all crossings of $pr_{wp}(D)$ have the winding parity $0$. In particular, for an oriented diagram $D$ only with crossing of winding parity $0$, $pr_{wp}(D)=D$. We obtain the following corollary.
\begin{cor}
    $pr_{wp} \circ pr_{wp} = pr_{wp}$.
\end{cor}

Now let us define a map $proj$ from $\mathcal{D}_{dl}^{0}$ to $\mathcal{VD}$ by removing all double lines. It is easy to see that two diagrams $D$ and $D'$ are equivalent up to the crossing sliding move and the moves in Fig. \ref{fig:moves2} except for the move (4), then $proj(D)$ and $proj(D')$ are equivalent up to generalized Reidemeister moves.
By Theorem~\ref{thm:wp-projection} and the map $proj$ we obtain the following statement.
\begin{cor}
    Let $D_{1}$ and $D_{2}$ be equivalent virtual knot diagrams as knots in $S_{g}\times S^{1}$. Then they are equivalent as virtual knots, that is, up to generalized Reidemeister moves.
\end{cor}

\begin{proof}
    Let $D_{1}$ and $D_{2}$ be equivalent virtual knot diagrams as knots in $S_{g}\times S^{1}$. Then $proj 
    \circ pr_{wp}(D_{i})= D_{i}$ for $i=1,2$. From Theorem~\ref{thm:wp-projection} and the previous observation on the images of $proj$, it follows that $proj 
    \circ pr_{wp}(D_{1})$ and $proj 
    \circ pr_{wp}(D_{2})$ are equivalent up to generalized Reidemeiter moves. Since $proj 
    \circ pr_{wp}(D_{i}) =D_{i}$, the proof is completed.
\end{proof}

Since classical knot theory is embedded well into the virtual knot theory, it follows that:
\begin{cor}
    Let $D_{1}$ and $D_{2}$ be equivalent classical knot diagrams as knots in $S_{g}\times S^{1}$. Then they are equivalent as classical knots, that is, up to Reidemeister moves.
\end{cor}

On the other hand, diagrams with double lines, which are equivalent to virtual knot diagrams, can be specified with the following properties.

\begin{thm}\label{thm:remove-doublelines}
       Let $D$ be an oriented diagram with double lines. Assume that $deg(D)=0$ and all crossings of $D$ have $0$ or $-1$ winding parity. Then $D$ is equivalent to $proj(D)$, that is, we can remove all double lines, by using only the crossing change move and the double line canceling move.
\end{thm}

\begin{proof}
    Let $D$ be an oriented diagram with double lines of $deg(D)=0$ such that all crossings of $D$ have $0$ or $-1$ winding parity. Then by crossing change moves, one can deform the diagram to a diagram with $0$ winding parity crossings. Without loss of generality, we assume that $D$ only has $0$ winding parity crossings.
    The remaining part of the proof consists of two steps.
    
    {\bf Step 1.} We will show that for a crossing $c$ one can remove all double lines on $\gamma_{c}$ by using crossing sliding move. For a crossing $c$, $\gamma_{c}$ consists of arcs. Let $\alpha_{1},\dots, \alpha_{k}$ be the sequence of arcs, which $\gamma_{c}$ consists of, such that $\alpha_{1}$ and $\alpha_{k}$ are arcs connected to $c$ and $\alpha_{i}$ and $\alpha_{i+1}$ cross a crossing $c_{i}$ on $\gamma_{c}$. Let the sum of signs of double lines on $\alpha_{i}$ be $d_{i}$ for $i=1,\dots, k$. First we apply crossing sliding moves on $c_{1}$ $|d_{1}|$ times to add double lines on $\alpha_{1}$ and $\alpha_{2}$ so that $\alpha_{1}$ and $\alpha_{2}$ have $-d_{1}$ and $d_{1}$ respectively. Then in the obtained diagram, there are double lines on $\alpha_{2}$ with $d_{1} + d_{2}$, all doubles on $\alpha_{1}$ can be canceled. Note that the remaining incoming and outgoing arcs of $c_{1}$ have $-d_{1}$ and $d_{1}$ double lines. Now we apply $|d_{1} + d_{2}|$ times crossing sliding move on $c_{2}$ so that we add $-(d_{1} + d_{2})$ double lines on $\alpha_{2}$. Then in the obtained diagram, $\alpha_{2}$ has no double lines and $\alpha_{3}$ has $d_{1} + d_{2}+d_{3}$ double lines. Note that the remaining incoming and outgoing arcs of $c_{2}$ have $-(d_{1} + d_{2})$ and $d_{1} + d_{2}$ double lines respectively.\\
    We repeat this process untill we meet one crossing twice. Let us say $c_{s} = c_{t}$ for $1 \leq s < t \leq k$. Note that the crossing $c_{t}$ consists of $\alpha_{s}$, $\alpha_{s+1}$, $\alpha_{t}$ and $\alpha_{t+1}$. By the previous process, $\alpha_{s}$ and $\alpha_{s+1}$ have no double lines, but $\alpha_{t}$ has double lines with $\Sigma_{i=s+1}^{t}d_{i}$. Since all crossings of $D$ have $0$ winding parity and $\alpha_{s+1} \cup \dots \cup \alpha_{t}$ is $\gamma_{c_{s}}$ or $K - \gamma_{c_{s}}$, $\Sigma_{i=s+1}^{t}d_{i}=0$. That is, $\alpha_{t}$ already has no double lines and so do $\alpha_s$ and $\alpha_{s+1}$. So we do not need to apply crossing sliding move on $c_{s}=c_{t}$. When we return to $c$, in the obtained diagram $\alpha_{k}$ has $\Sigma_{i=1}^{k}d_{i}$ double lines and it is winding parity of $c$, which is $0$.
    
    {\bf Step 2.} Let us take any crossing $c$. 
    From Step 1., we can remove all double lines on $\gamma_{c}$. Now we follows $\gamma_{c}' = D-\gamma_{c}$ starting from $c$. Assume that we meet a crossing $c'$ on $\gamma_{c}'$. If it appeared in $\gamma_{c}$, then we can see that the arc before $c'$ has no double lines, since winding parities of all crossings are $0$. So we pass $c'$ without any deformation. Suppose that $c'$ is not contained in $\gamma_{c}$. Then we apply the process in Step 1. Note that after this process, no new double lines on $\gamma_{c}$ appear. We continue the process; we do the process in Step 1, when we meet a crossing, which is not appeared in $\gamma_{c}$, but we pass crossings, which appeared in $\gamma_{c}$ already.\\
    When we meet $c'$ again on $\gamma_{c}'$, the arc incoming to $c'$ has double lines, the number of which is the winding parity of $c'$ by the argument in Step 1. 
    By repeating this process, we can show that it is possible to remove all double lines by crossing sliding moves and by canceling double lines. It completes the proof.
\end{proof}

\begin{cor}
    Let $K$ be a knot in $S_{g} \times S^{1}$ of degree $0$. Assume that $K$ is placed in $\Sigma_{1} \times S^{1} \subset S_{g} \times S^{1}$ such that $S_{g} = \Sigma_{1} \cup \Sigma_{2}$ and its diagram with double lines have all crossings with $0$ or $-1$ winding parity. Then there exists an ambient isotopy $\{h_{t}\}_{t\in [0,1]}$ such that $h_{t} | \Sigma_{2} \times S^{1} = id_{\Sigma_{2} \times S^{1}}$ and $h_{1}(K) \in \Sigma_{1} \times [0-\epsilon, 0+\epsilon]$ for some $\epsilon>0$.
\end{cor}

\begin{proof}
    Suppose that $K$ is a knot of degree $0$ placed in $\Sigma_{1} \times S^{1}  \subset S_{g} \times S^{1}$ and a diagram with double lines of $K$ have all crossings with $0$ or $-1$ winding parity. Then by Theorem~\ref{thm:remove-doublelines} we can obtain a diagram without double lines by only the crossing change move, the crossing sliding move and the double line canceling. Note that each move can be presented by isotopy inside $\mathbb{D}^{2} \times S^{1} \subset \Sigma_{1} \times S^{1}  \subset S_{g} \times S^{1}$ such that $\mathbb{D}^{2}$ is a disc, where these moves are applied. It completes the proof.
\end{proof}

\section{Essential diagrams for links in $S_{g} \times S^{1}$}

Let us denote an oriented diagrams with double lines $\dL$ with a base diagram $D$ by $D_{\dL}$. A subset $\dd \subset \dL$ is called {\em a set of important double lines} if $D_{\dL - \dd}$ is a diagram of degree $0$ such that all crossings of it have degree $0$ or $-1$. If a set $\dd$ of important double lines has the minimal number of double lines among all possible sets of important double lines, then it is called {\em a set of essential double lines.} 

\begin{exa}
    Let $D_{\dL}$ be a diagram with double lines. Then $\dL$ is a set of important double lines. 
\end{exa}

\begin{figure}
    \centering
    \includegraphics[width = 2cm]{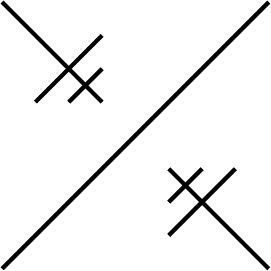}
    \caption{}
    \label{fig:minus-important}
\end{figure}

\begin{lem}
    Let $D_{\dL}$ be a diagram with double lines. Let $\dd$ be a set of essential double lines. Then $D_{\dL}$ can be deformed to a diagram $D_{\dL'}$ so that $\dL' = \dd$ or $\dL' = \dd \cup \dd_{-1}$, where double lines in $\dd_{-1}$ appear near to crossings with winding parity $-1$ in $D_{\dL - \dd}$ as described in Fig.~\ref{fig:minus-important}. We call such a diagram $D_{\dL'}$ {\em an essential diagram}.
\end{lem}

\begin{proof}
    Let $D_{\dL}$ be a diagram with double lines. Let $\dd$ be a set of essential double lines. By Theorem~\ref{thm:remove-doublelines}, $D_{\dL - \dd}$ can be deformed to $D'$ without double lines. If $D_{\dL - \dd}$ has crossings only with winding parity $0$, then $D'=D$. If not, $D'$ can be deformed to $D_{\dd_{-1}}$ by applying crossing change moves on crossings with winding parity $-1$ in $D_{\dL - \dd}$. Since we only use crossing change moves, crossing sliding moves and double line canceling in this process, it can be done for $D_{\dL}$ to obtain $D_{\dd\cup\dd_{-1}}$ and it completes the proof.
\end{proof}
\begin{rem}
    Note that $D_{\dL}$ and $D_{\dL'}$ have the same base diagram.
\end{rem}

\begin{rem}
    In the notion of essential diagrams, a diagram with double lines of degree $0$ such that all crossings have the winding parity $0$, has the empty set as a set of essential double lines and virtual diagrams as the essential diagram.
\end{rem}

\begin{lem}\label{lem:nbr-ess-dl}
    Let $D_{\dd}$ be an essential diagram. Assume that $D_{\dd}$ is equivalent to $D'_{\dL'}$. Then the number of essential double lines of $D'_{\dL'}$ is less than or equal to $|\dd|$.
\end{lem}
\begin{proof}
 Assume that $$D_{\dd}=D_{0} \rightarrow D_{1} \rightarrow \dots \rightarrow D_{n-1} \rightarrow D_{n} = D'_{\dd'},$$
 where $D_{i+1}$ is obtained from $D_{i}$ by one of moves. For each $D_{i+1}$ we give a set $\dL_{i+1}$ of important double lines induced from the given set $\dL_{i}$ of important double lines of $D_{i}$. 

 If $D_{i+1}$ is obtained from $D_{i}$ by one of generalized Reidemeister moves, then we take $\dL_{i+1} = \dL_{i}$ for $D_{i+1}$. It is clear that $\dL_{i+1}$ is still a set of important double lines.

 If $D_{i+1}$ is obtained from $D_{i}$ by the move (4), we take $\dL_{i+1} = \dL_{i}$ for $D_{i+1}$. If the double line in the move $(4)$ is an important double line, then $\dL_{i+1}$ is a set of important double lines. If the double line in the move $(4)$ is not an important double line, then the winding parity of the crossing in the move (4) is changed. But, if the crossing has a winding parity $0$ in $D_{i}$ after removing $\dL_{i}$, then the crossing has a winding parity $-1$ in $D_{i+1}$ after removing $\dL_{i+1}$. Similarly, if the crossing has a winding parity $-1$ in $D_{i}$ after removing $\dL_{i}$, then the crossing has a winding parity $0$ in $D_{i+1}$ after removing $\dL_{i+1}$. Therefore, $\dL_{i+1}$ is still a set of important double lines.
 
 If $D_{i+1}$ is obtained from $D_{i}$ by appearing two neighboring double lines, then take $\dL_{i+1} = \dL_{i}$. It is clear that $\dL_{i+1}$ is still a set of important double lines.

 Assume that $D_{i+1}$ is obtained from $D_{i}$ by double line canceling. If both of them are not important in $\dL_{i}$, then take $\dL_{i+1} = \dL_{i}$. If both of them $d,d'$ are in $D_{i}$, then take $\dL_{i+1} = \dL_{i} - \{d,d'\}$. One can show that $\dL_{i+1}$ is a set of important double lines in $D_{i+1}$ in both cases. Now suppose that $d$ is important, but $d'$ is not important. Note that if $d'$ is not important, then it was in $\dd_{-1}$ of $D_{\dd}$ or it appeared by double line appearing move. In both cases, there is a naturally paired double line $d''$ to $d'$ with a different sign with $d'$, but the same sign with $d$. Now take $\dL_{i+1} = (\dL_{i} - \{d\}) \cup \{d''\}$. It is clear that the sum of signs of double lines in $\dL_{i+1}$ is $deg(D_{\dL})$. 

\textbf{Claim.} $\dL_{i+1}$ is a set of important double lines: Let us $D_{s}$, $0\leq s\leq i$ be a diagram such that $d'$ and $d''$ appear. Following our construction, $D_{s}$ has a set $\dL_{s}$ of important double lines, but $d',d'' \not\in \dL_{s}$. Let us take $\dL_{s}' = \dL_{s} \cup \{d',d''\}$. Then $\dL_{s}'$ is a set of important double lines again.  Moreover, it induces another set $\dL_{i}'$ of important double lines for $D_{i}$ such that $\dL_{i}' =dL_{i} \cup \{d',d''\}$. Now $d, d',d''$ are important double lines and our situation becomes that two important double lines are canceled. Then $\dL_{i+1} = \dL_{i}'- \{d,d'\}$ and it is a set of important double lines of $D_{i+1}$.

Now all diagrams $D_{i}$ have corresponding sets of important double lines. Note that $|dL_{i}| \geq |dL_{i+1}|$ for all $i = 0, \dots, n$. Hence the number of essential double lines of $D'_{\dL'}$ is less than or equal to $|\dL_{n}| \leq |\dd|$. 
\end{proof}

\begin{thm}
    Let $D_{\dd}$ and $D_{\dd'}'$ be equivalent essential diagrams. Then the number of essential double lines of $D_{\dd}$ and $D_{\dd'}'$ are the same.
\end{thm}

\begin{proof}
    It follows from Lemma~\ref{lem:nbr-ess-dl}.
\end{proof}


\section{Application to study classical links}
In this section, we discuss several applications of essential diagrams to study classical links of 2-components. 

\subsection{Knots with one crossing of degree $k \neq 0$}

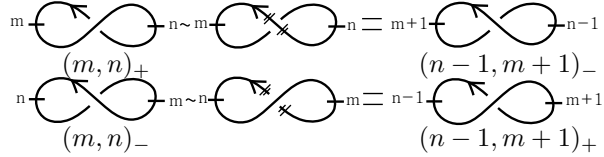
\begin{figure}[h]
\begin{center}
\includegraphics[width = 8cm]{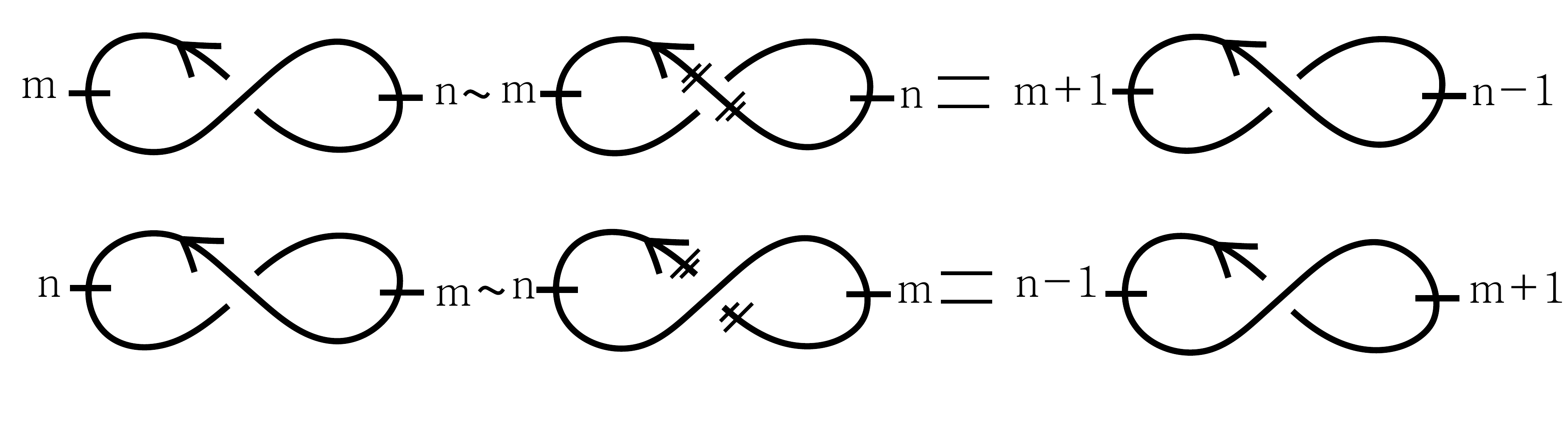}
\put(-205, 30){$(m,n)_{+}$}
\put(-205, 2){$(m,n)_{-}$}
\put(-70, 30){$(n-1,m+1)_{-}$}
\put(-70, 2){$(n-1,m+1)_{+}$}
\end{center}
\caption{Knot diagrams with one crossing corresponding to $(m,n)_{\epsilon}$ and their equivalence relation}\label{fig:m_n-one-crossing}
\end{figure}

For integers $m,n$, let us denote a diagram with one crossing of crossing sign $\epsilon$ by $(m,n)_{\epsilon}$ if the arc from under-crossing to over-crossing has a line with an integer number $m$ and remaining part has a line with an integer number $n$, as described in Fig. \ref{fig:m_n-one-crossing}. The degree of $(m,n)_{\epsilon}$ is equal to $m+n$.

From the moves for diagrams with double lines, we obtain that $(m,n)_{\epsilon}$ is equivalent to $(n-1,m+1)_{-\epsilon}$. In \cite{Kim-winding-parity} we partially classified $(m,n)_{\epsilon}$ for $m,n \in \mathbb{Z}$.

\begin{lem}\label{lem:(m,n)-essential}
    Let $m,n \in \mathbb{Z}$. Then the diagram $(m,n)_{\epsilon}$ is an essential diagram.
\end{lem}

\begin{proof}
    For $m,n \in \mathbb{Z}$ $(m,n)_{\epsilon}$ has two halves $\gamma_{c}$ and $\gamma_{c}'$ at $c$ such that $\gamma_{c}$ starts from under-crossing of $c$ and ends at the over-crossing of $c$, but $\gamma_{c}'$ is remaining arc. Note that double lines on $\gamma_{c}$ determines the winding parity of $c$. Let $\dL$ be the set of all double lines on $(m,n)_{\epsilon}$ and $\dd$ a set of important double lines. Let us denote $D_{\dL} = (m,n)_{\epsilon}$. 

    {\bf Case 1.} $m \geq 0$.\\
    Since $D_{\dL-\dd}$ has one crossing with winding parity $0$ or $-1$, there should be no double lines or only one double line with sign $-1$ on $\gamma_{c}$ for $D_{\dL-\dd}$. Then all double lines on $\gamma_{c}$ for $D_{\dL}$ are contained in $\dd$ or all but one double lines with sign $-1$ on $\gamma_{c}$ for $D_{\dL}$ belong to $\dd$. Since $m>0$, there is no double line with sign $-1$ on $\gamma_{c}$ for $D_{\dL}$, and hence all double lines on $\gamma_{c}$ for $D_{\dL}$ belong to $\dd$. Moreover, since the sum of signs of double lines in $\dd$ is equal to $m+n$, all double lines on $\gamma_{c}'$ also should be contained in $\dd$. It follows that $\dd=\dL$. That is, $\dd$ is the set of essential double lines and $(m,n)_{\epsilon}$ is an essential diagram.

    {\bf Case 2.} $m \leq -1$.\\
    Similarly, all double lines on $\gamma_{c}$ for $D_{\dL}$ are contained in $\dd$ or all but one double lines with sign $-1$ on $\gamma_{c}$ for $D_{\dL}$ belong to $\dd$. If all double lines on $\gamma_{c}$ are contained in $\dd$, then all double lines on $\gamma_{c}'$ also should be in $\dd$. Then the number of double lines is $-m+|n|$. Let us denote this set of important double lines by $\dd^{0}$. Let us denote this set of important double lines by $\dd^{1}$.
    
    Assume that all double lines but one on $\gamma_{c}$ are contained in $\dd$. Since the sum of signs of double lines in $\dd$ is $m+n$, all but one double lines on $\gamma_{c}'$ with sign $+1$ should be contained in $\dd$. Then $\dd$ has $|m|-1 + |n|-1 = -m+|n|-2$ and the sum of signs of double lines is $m+1 +n -1 = m+n$. 
    
    If $n\leq 0$, there does not exist $\dd^{1}$, that is, $\dd^{0} =\dL$ is one and only possible set of important double lines. Therefore $\dL$ is the set of essential double lines and $(m,n)_{\epsilon}$ is an essential diagram.
    
    If $n>0$, then there are two possible sets $\dd^{0}$ and $\dd^{1}$ of important double lines. Since $|\dd^{1}| < |\dd^{0}|$, $\dd^{1}$ is the set of essential double lines. Moreover, $D_{\dL-\dd^{1}}$ has two double lines placed as described in Fig.~\ref{fig:minus-important}. $(m,n)_{\epsilon}$ is an essential diagram.
    \end{proof}

As a corollary, one can obtain the following statement.
\begin{cor}
    Let $m,n,m',n'$ be integers. Let $\dd((m,n)_{\epsilon})$ be the number of essential double lines. If $\dd((m,n)_{\epsilon}) \neq \dd((m',n')_{\epsilon'})$, then $(m,n)_{\epsilon}$ and  $(m',n')_{\epsilon}$ are not equivalent.
\end{cor}

From the above corollary, we can expand the partial classification of $(m,n)_{\epsilon}$. In \cite{Kim-winding-parity} the following statement is proved by using winding parity:
\begin{prop}\label{prop:infinite-nontrivial-1cro-degk}
For an integer $k$ with $k\geq 3$, there are at least
  \begin{equation*}
    \left\{ \begin{array}{lr}
      \lfloor \frac{k-2}{2} \rfloor +1 & \text{if}~ $k$ ~ \text{is odd} \\
       \lfloor \frac{k-2}{2} \rfloor       & \text{if}~ $k$ ~ \text{is even}
    \end{array}\right.
  \end{equation*}
nontrivial knots of degree $k$ with one crossing.
\end{prop}

But two knots $(2,k-2)_{+}$ and $(2+mk,k-mk-2)_{+}$ cannot be distinguished by using winding parity for $k\geq3$. Now by using essential double lines, one can prove the following statement.

\begin{cor}
    Let $m,k$ be integers with $k\geq 3$ and $m\not\equiv 0$ mod $k$. For any positive integer $s$, $(m,k-m)_{\epsilon}$ and $(m+sk,k-sk-m)_{\epsilon}$ are not equivalent. That is, for each $k$, there are infinitely many nontrivial knots of degree $k\geq 3$ with one crossing.
\end{cor}

\begin{proof}
    By the proof of Lemma~\ref{lem:(m,n)-essential}, $(m,k-m)_{\epsilon}$ has the set $\dd$ of essential double lines such that $|\dd|= |m| + |k-m|$. If $m+sk\geq 0$, then the set $\dd'$ of essential double lines of $(m+sk,k-sk-m)_{\epsilon}$ has $|m+sk| + |k-sk-m|$. Since $|\dd'|= m+sk + |k-sk-m|\geq m+sk+|k-m|-|sk| = m + |k-m| = |\dd|$. Since $m \not\equiv 0$ mod $k$, we obtain $|\dd'| > |\dd|$ and it completes the proof.
\end{proof}

\subsection{2-component links with the trivial component and knots in $S^{2} \times S^{1}$}
We denote by $T$ a trivial knot. For a two component link $L = K \sqcup T$, it is clear that $K \in S^{3}-N(T)$ is a knot in $D^{2} \times S^{1}$. In particular, $K$ can be considered as a knot in $S^{2} \times S^{1}$, which is the result of 2-handle surgery along $T$ with $0$-framing. It is clear that, if there exists an ambient isotopy between two links $L$ and $L'=K' \sqcup T$, then the corresponding knots $K$ and $K'$ in $S^{2}\times S^{1}$ are ambient isotopic, that is, equivalent. In this section, we study 2-component links $L = K \sqcup T$ with a trivial knot $T$ by using knots in $S^{2} \times S^{1}$.

\begin{figure}
    \centering
    \includegraphics[width=0.5\linewidth]{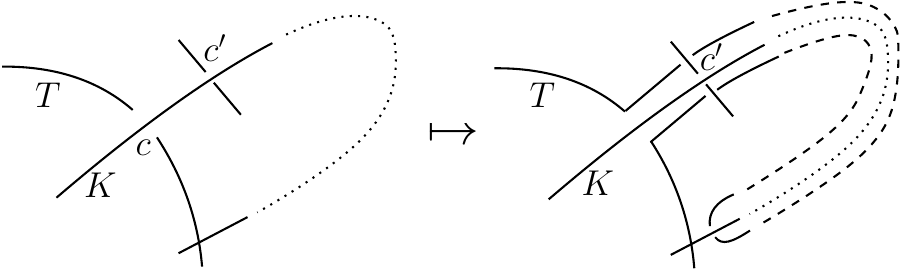}
    \caption{We deform the trivial component so that all crossings are placed inside the disc, which is bounded by the trivial component}
    \label{fig:exa-deform-to-dldiag}
\end{figure}

Let us describe a way to obtain a diagram with double lines from $L = K \sqcup T$. Let $D$ be a diagram of $L$ such that the component of $D$ corresponding to $T$ is a trivial circle. With an abuse of the notation, we denote the component of $D$ corresponding to $K$ and $T$ by $K$ and $T$, respectively. Assume that the trivial component $T$ is oriented in counter-clockwise orientation. We call the infinite region of $\mathbb{R}^{3} - T$ to be {\em outside of $T$.} The remaining region is called to be {\em inside of $T$.} One can obtain a diagram with double lines for $K$ from the diagram $D$ as a knot in $D^{2} \times S^{1}$ as follows:
Take a crossing $c$ between $K$ and $T$. Suppose that an arc of $T$ is the under-crossing of $c$. Then the over-crossing $\beta$ of $c$ passes $T$ from inside to outside of $T$. We push the under-crossing $\alpha$ of $c$ along the over-crossing arc to outside. If we meet a crossing $c'$ between $K$ and itself, then we push $\alpha$ so that it passes under the other three arcs connected with $c'$, see Fig.~\ref{fig:exa-deform-to-dldiag}. 
We continue this process till we meet another crossing between $K$ and $T$. This process can be done when $T$ meets a crossing $c$ between $K$ and $T$ as an over-crossing.

If there is no crossing outside $T$, then we stop our process. If there is a crossing $c''$ between $K$ and itself outside $T$, then we take another crossing between $K$ and $T$ closest to $c''$ and repeat the above process. Since the number of crossings is finite, this process ends in a finite step. Note that we do not make new crossings between $K$ and itself and the linking number is not changed.

\begin{figure}
    \centering
    \includegraphics[width=0.5\linewidth]{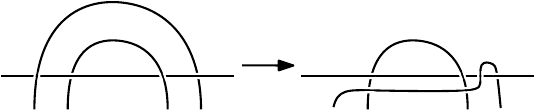}
    \caption{We deform an arc of $K$ connecting two crossings between $K$ and $T$ so that all pairs of crossings between $K$ and $T$ connected by arcs of $K$ are neighboring.}
    \label{fig:deformation-sewed}
\end{figure}

Now all crossings $K$ are placed inside $T$ and there are finitely many arcs of $K$ connecting two crossings between $K$ and $T$. Take one arc. If it connects two crossings so that it is under-crossing (or over-crossing) of both crossings, then we push it inside $T$ by using the second Reidemeister move. If it connects two crossings alternatively, then we push the arc inside $T$ so that two end are neighboring, for example, see Fig.~\ref{fig:deformation-sewed}.
Since there are finitely many arcs outside, we can deform the diagram so that all outside arc connect crossings alternatively. We say that {\em $T$ is sewed with $K$ in this diagram.}

\begin{figure}
    \centering
    \includegraphics[width=0.5\linewidth]{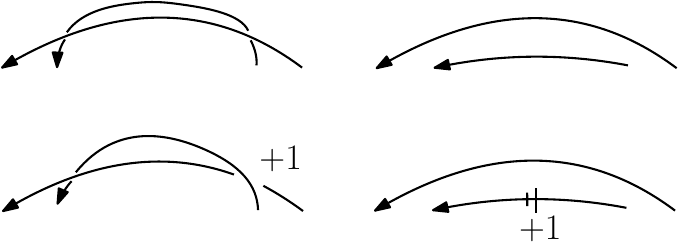}
    \put(-40,40){-1}
    \put(-112, 56){-1}
    \caption{We replace a linked part to a double line with respect to the linking number}
    \label{fig:change-to-dldiagram}
\end{figure}

Now we replace alternatively connecting arc by an arc inside $T$ with double line such that it has the linking number coming from those two crossings as the sign of a double line, see Fig.~\ref{fig:change-to-dldiagram}. 
Finally we obtain a diagram with double line on the disc bounded by $T$. Let us denote it by $D_{K}^{dl}$. Since we consider 2-component links with trivial component, which is placed really well, it is clear that this map is well defined.

\begin{rem}
    Note that one difference between knots in $D^{2} \times S^{1}$ and $S^{2} \times S^{1}$ is the possibility of ``the infinite point change'' move. In the sense of links $L = K \sqcup T$ in $S^{3}$ it can be presented by the handle sliding move along $T$ with $0$-framing.
\end{rem}


By using the above correspondence, now we can use knots in $S_{g}\times S^{1}$ to study a family of classical links.

\subsection{2-component links with one crossing}
\begin{figure}
    \centering
    \includegraphics[width=0.5\linewidth]{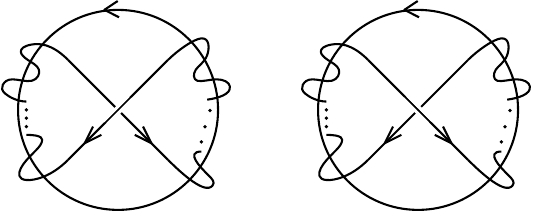}
    \put(-190, 30){$m$}
    \put(-102, 30){$n$}
    \put(-83, 30){$m$}
    \put(0, 30){$n$}
    \put(-148, -10){$L_{(m,n)_{+}}$}
    \put(-46, -10){$L_{(m,n)_{-}}$}
    \caption{$L_{(m,n)_{\epsilon}}$ has a diagram of 2-component links such that the self crossing has crossing sign $\epsilon$ and linking numbers of left and right part of the diagram are $m$ and $n$ for $m,n \in \mathbb{Z}$}
    \label{fig:Lmn}
\end{figure}
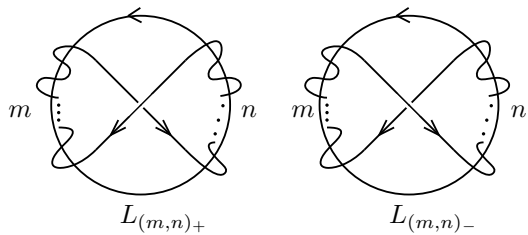

Let us denote the link described in Fig.~\ref{fig:Lmn} by $L_{(m,n)_{\epsilon}}$. It is clear that $L_{(m,-m)_{\epsilon}}$ has the linking number $0$.

\begin{cor}
    For any $m,n \geq 1$ with $m\not=n$, $L_{(m,-m)_{+}}$ and $L_{(n,-n)_{+}}$ are not equivalent. Moreover, the number of crossings is minimal. 
\end{cor}

\begin{proof}
    From $L_{(m,-m)_{+}}$ one can obtain a diagram with double lines as described in Fig.~\ref{fig:m_n-one-crossing}, denoted by $(m,-m)_{+}$. In \cite{Kim-winding-parity} it is shown that $(m,-m)_{+}$ is not trivial. Moreover, $(m,-m)_{+}$ is equivalent to $(-m+1,m-1)_{-}$, but the number of double lines and crossing is not changed and hence, it completes the proof.
\end{proof}

\begin{thm}
    Let $L = K \sqcup T$ be a 2-components link such that $T$ is a trivial knot. Assume that $lk(K,T)=0$ and $T$ is sewed by $K$. If $lk(\gamma_{c},T)=0, -1$ for all crossings $c$, where $\gamma_{c}$ is the half at the crossing $c$, then $K$ and $T$ can be separated.
\end{thm}

\begin{proof}
    Let $D_{K}^{dl}$ be the diagram with double lines obtained from the above process. It is easy to see that $deg(D_{K}^{dl})=0$. Moreover, there is natural correspondence between crossings of $K$ and $D_{K}^{dl}$. From the construction, we obtain that the winding parity of $c$ for each crossing of $D_{K}^{dl}$ is $0$ or $-1$. By Theorem~\ref{thm:remove-doublelines}, $D_{K}^{dl}$ can be deformed to a diagram without double lines by using crossing sliding move. Note that the crossing sliding move presents an ambient isotopy $h_{t}$ of $S^{2} \times S^{1} = \mathbb{D}_{1}^{2} \times S^{1} \cup_{id|S^{1} \times S^{1}} \mathbb{D}_{2}^{2} \times S^{1}$ such that $h_{t}|_{\mathbb{D}_{2}^{2}} = id_{\mathbb{D}_{2}^{2}}$. That is, we obtain an ambient isotopy of $S^{3}-N(T)$ such that it place $K$ seperated from $T$. It completes the proof.
\end{proof}

\begin{rem}
    In this process, the diagram of $T$ is placed in a good position. So, `$T$ is sewed by $K$' means that the diagram of $K$ is placed on the disc of $T$ well except finitely many arcs linked with $T$. 
\end{rem}

\end{document}